\documentclass{amsart}
\usepackage{amssymb,amscd}
\title{Deformation of Weyl Modules\\ and Generalized Parking Functions}
\author{B. Feigin, S. Loktev}
\address{BF: Landau institute for Theoretical Physics, Chernogolovka
142432, Russia} \email{feigin@feigin.mccme.ru}
\address{SL:  Department of Mathematics, Kyoto University, Kyoto 606-8502, Japan}
\address{On leave from Institute for Theoretical and Experimental Physics;
Independent University of Moscow}
\email{loktev@math.kyoto-u.ac.jp}
\date{}
\newtheorem{prp}{Proposition}
\newtheorem{lmm}{Lemma}
\newtheorem{thm}{Theorem}
\newtheorem{cnj}{Conjecture}
\newtheorem{crl}{Corollary}

\newcommand{\bnr}[2]{\{1,\dots,#1\}\times \{1, \dots, #2\}}
\begin{document}
\def\theenumi{\roman{enumi}}
\def\labelenumi{(\theenumi)}

\begin{abstract}
Local Weyl modules over two-dimensional currents with values in $gl_r$
are deformed into spaces with bases 
related to parking functions. Using this construction we

1) propose a simple proof that dimension of the space of diagonal coinvariants is not less than
the number of parking functions;

2) describe the limits of Weyl modules in terms of semi-infinite forms and find the
limits of characters;

3) propose a lower bound and state a conjecture for
dimensions of Weyl modules with arbitrary highest weights.

Also we express characters of deformed Weyl modules
in terms of $\rho$-parking functions and the Frobenius characteristic map.
\end{abstract}

\maketitle

\def \CC {{\mathbb C}}
\def \ZZ {{\mathbb Z}}
\def \NN {{\mathbb N}}
\def \g {{\mathfrak g}}
\def \h {{\mathfrak h}}
\def \n {{\mathfrak n}}
\def \b {{\mathfrak b}}
\def \ad {{\rm ad}\,}
\def \Z {{\mathcal Z}}

\section*{Introduction}

\def \bo {\mathfrak b}
\def \al {\mathfrak a}

\def \slt {\mathfrak{sl}_2}
\def \glt {\mathfrak{gl}_2}
\def \slth {\widehat{\slt}}
\def \hom {{\rm {Hom}}}
\def \gl {{\mathfrak gl}}
\def \hk {\hookrightarrow}
\def \V {{\mathcal V}}

In this paper we continue studying Weyl modules in the sense of \cite{CP}, \cite{FL2}.

Fix a commutative algebra $A$ with unit, an augmentation $\theta: A \to \CC$
and a semi-simple Lie algebra $\g$ with a Cartan decomposition
$\g=\n_- \oplus\h\oplus\n_+$. Let $\bo=\n_+ \oplus \h$.

Recall that local Weyl modules are cyclic representations of the
Lie algebra $\g\otimes A$ labeled by dominant characters
$\chi: \h\to \CC$.
Namely  they are the quotients $I_\chi/J_\chi$, where
$I_\chi$ is the $\g \otimes A$-module induced from the following $1$-dimensional representation
of $\bo\otimes A$:
$$\bo\otimes A\stackrel{1\otimes \theta}{\longrightarrow} \bo\to\h
\stackrel{\chi}{\to} \CC.$$
and $J_\chi$ is the minimal submodule in $I_\chi$ such that $I_\chi/J_\chi$
is integrable as a representation of $\g\simeq \g\otimes 1\subset \g \otimes A$.
For a generic $A$ these modules are studied in \cite{FL2}, but still
almost nothing is known about their structure.

\medskip

First let us describe what is known in the simplest case of
$A=\CC[x]$, the augmentation $\theta: P(x) \to P(0)$ and $\g=\slt$. 
Then the Weyl modules are labeled by a non-negative integer 
and denoted in \cite{FL2}  by $W_\CC (\{0\}_{n\omega})$. For a shortness let us denote them
here by $W_1(n)$.
It is shown in \cite{CP} that $\dim W_1(n)=2^n$.
Let us outline a couple of ways to describe these modules.

\medskip

\noindent {\bf First way.} Such a
construction is introduced  in \cite{FL1}
and is called {\em fusion}.

For $z \in \CC$ 
let $\CC^2(z)$ be the evaluation two-dimensional representation of
$\slt\otimes \CC[x]$. Recall that an element $R\in\slt\otimes \CC[x]$ 
acts on $\CC^2(z)$ just by the $2\times 2$-matrix $R(z)$.
By $v(z)$ denote the highest
weight vector of $\CC^2(z)$.

 Consider the tensor product
$$W(z_1,\ldots, z_m) = \CC^2(z_1)\otimes\ldots\otimes \CC^2(z_m).$$
The module $W(z_1,\ldots,z_m)$ is cyclic (even irreducible) if $z_i\ne z_j$ for $i\ne j$.
Choose $v(z_1)\otimes \ldots \otimes v(z_m)$ as a cyclic vector.
Then the family $W(z_1,\ldots, z_m)$ has a limit when all
$z_j\to 0$ in the class of cyclic modules over $\slt\otimes \CC[x]$ (see \cite{FF} for the proof). 
Combining the results of \cite{FF} and \cite{CP} one can see that
this limit is isomorphic to $W_1(m)$. 

\medskip

\def \cl {{\mathcal L}}
\noindent {\bf Second way.} 
Let $\slth$ be the central extension of the Lie algebra $\slt\otimes \CC [x,x^{-1}]$.
Then $W_1(m)$ can be realized as a {\em Demazure
submodule } of a fundamental integrable representation of $\slth$.

By $\cl_0$ and $\cl_1$ denote the fundamental integrable representations
of the affine Lie algebra $\slth$.
Let $h_0 = h \otimes 1 \in \slth$, where $h\in \slt$ is the generator  of the
Cartan subalgebra. Introduce the notation $v_m$ for the extremal vectors of $\cl_0$ and
$\cl_1$ such that $h_0 v_m=m v_m$, so $v_m \in \cl_0$ for even $m$ and
$v_m \in \cl_1$ for odd $m$. Combining the results of \cite{FL1} and \cite{CP}
one can see that the submodule generated by $v_m$ under the action of
the annihilating operators $\slt\otimes \CC[x]\subset \slth$ is isomorphic to $W_1(|m|)$.

In particular it
follows that $W_1(m)$ is dual to the space of sections of some line
bundle on the affine $\slth$-Schubert variety of dimension $m$.

Such a description of the Weyl modules demonstrates that the latter can be
used to construct the fundamental representations of $\slth$. 
Let $M$ be a $\slt\otimes \CC[x]$-module. By definition of Weyl module
homomorphisms $W_1(m)\to M$ are in
one-to-one correspondence with elements $v\in M$
such that 
$$(e\otimes x^i)\cdot v=0\ \ \mbox{for} \ i\ge 0; \qquad
 (h\otimes x^j)\cdot v=0 \ \ \mbox{for} \ j> 0; \qquad (h\otimes 1)\cdot v=mv.$$ 
Using this property one can show that
$\hom (W_1(m), W_1(m+2))$ is one-dimensional and spanned by an embedding.
Writing down two sequences of homomorphisms
\begin{gather}
\label{0}
W_1(0)\hk W_1(2)\hk W_1(4)\hk\ldots,\\
\label{1}
W_1(1)\hk W_1(3)\hk W_1(5)\hk\ldots,
\end{gather}
one can show that $\cl_0$ and $\cl_1$ are the inductive limits of $\eqref{0}$ and 
$\eqref{1}$ respectively. Note that despite of each $W_1(m)$ has only the structure of
$\slt\otimes \CC[x]$-module, the limit has a more reach structure of $\slth$-module.

\medskip

In this paper we consider the case of $A=\CC[x,y]$ and again 
$\theta: P(x,y) \to P(0,0)$. The corresponding Weyl modules are denoted by 
$W_{\CC^2}(\{0\}_\lambda)$ 
in \cite{FL2}. Dimensions of some of these modules are calculated in \cite{FL2} using
M.~Haiman's results about two-dimensional diagonal coinvariants \cite{H}.

Here let us restrict ourselves to the case $\g = \slt$, so 
Weyl modules are also labeled by an integer. Denote them by $W_2(m)$. It is shown in
\cite{FL2} that their dimensions are equal to Catalan numbers.

Concerning the first way of description, note that the algebra
$\slt\otimes \CC[x,y]$ also has the evaluation representations $\CC^2(Z)$,
where $Z$ is a point on the complex plane. The tensor product
$\CC^2(Z_1)\otimes\ldots\otimes \CC^2(Z_m)$ is cyclic 
for pairwise distinct $Z_i$ and we can consider the similar
limit when $\{Z_i\to (0,0)\}$. But in this case the limit depends on how the points
$\{Z_j\}$ go to the origin. Namely the set of all limits can be identified
with the special fiber $H_m^{0}$ in the Hilbert scheme of $m$-points on the plane.
Moreover the limits form a $2^m$-dimensional vector bundle $\mu_m$ over $H_m^0$ with the structure
of $\slt\otimes \CC[x,y]$-module on each fiber. It can be obtained
from the $m!$-dimensional bundle constructed in \cite{H} by the fiber-wise
Schur-Weyl duality.
Combining the results of \cite{FL2} and \cite{H} one can see that 
the Weyl module $W_2(m)$ is isomorphic to the space of sections of the bundle $\mu_m$.

This paper makes possible the second way of description.
Similar to the $1$-dimensional case these Weyl modules can be mapped
to each other. But in this case the space $\hom (W_2(m), W_2(m+2))$
is $m+2$-dimensional. Therefore there are many inductive limits like
$$W_2(0)\stackrel{\alpha_0}{\longrightarrow}
W_2(2)\stackrel{\alpha_2}{\longrightarrow}
W_2(4)\stackrel{\alpha_4}{\longrightarrow}\ldots \quad\text{ or }\quad
W_2(1)\stackrel{\alpha_1}{\longrightarrow}
W_2(3)\stackrel{\alpha_3}{\longrightarrow}\ldots.
$$
The limits are infinite-dimensional $\slt\otimes \CC[x,y]$-modules. In this paper we consider
rather special choice of the homomorphisms $\{\alpha_j\}$.
To describe it note that the $U(\slt\otimes \CC[x])$-submodule in $W_2(m)$ generated by
the cyclic vector is isomorphic to the Weyl module $W_1(m)$. 
And there is the unique way to choose the maps $\{\alpha_j\}$ such that the following diagrams 
are commutative:
$$
\begin{CD}
W_2(0) @>>> W_2(2) @>>> W_2(4) @>>> \ldots\  \cl_{0}^{(2)}\\
@AAA    @AAA  @AAA   @AAA\\
W_1(0) @>>> W_1(2) @>>> W_1(4) @>>> \ldots\  \cl_{0},
\end{CD}
$$
$$
\begin{CD}
W_2(1) @>>> W_2(3) @>>> W_2(5) @>>> \ldots\  \cl_{1}^{(2)}\\
@AAA    @AAA  @AAA   @AAA\\
W_1(1) @>>> W_1(3) @>>> W_1(5) @>>> \ldots\  \cl_{1}.
\end{CD}
$$
In this paper we show that the action of $\slt\otimes \CC[x,y]$ on the inductive limits 
$\cl^{(2)}_{0}$ and $\cl^{(2)}_{1}$  can be extended to an action of 
the universal central extension
$(\slt\otimes \CC[x,x^{-1},y])^{\widehat{}}$ 
of the algebra $\slt\otimes \CC[x,x^{-1},y]$. Let us briefly describe this algebra.

Let $\g$ be a simple Lie algebra and let $M$ be an affine variety. 
By $C(M)$ denote the ring of functions on $M$, by $\Omega^i$
denote the space of the $i$-forms on $M$.
Then we have
$H_2(\g\otimes C(M))\simeq \Omega^1/d\Omega^0$, so 
$$(\g\otimes C(M))^{\widehat{}} = \g\otimes C(M) \oplus \Omega^1/d\Omega^0$$
with the cocycle given by the formula
$$\omega(g_1\otimes P_1, g_2\otimes P_2) =
(g_1,g_2)\cdot \pi (d P_1\cdot P_2),$$
where $(,)$ is the Cartan-Killing form and  $\pi$ is the projection
$\Omega^1\to \Omega^1/d\Omega^0$. In our case $M = \CC^* \times \CC$.

\def \x {{\mathcal X}}
\def \y {{\mathcal Y}}

The main tool we use in this paper is deformation of the Weyl modules. After deformation
they become representations of the matrix algebra with values in differential operators 
on the line preserving order of zero at the origin.
More precisely, the polynomial ring $\CC[x,y]$ is deformed into the associative
algebra $\CC\left<\x,\y\right>$ generated  by the elements $\x$ and $\y$
(corresponding to
$x$ and $x\partial/\partial x$ in ${\rm Diff}(x)$ respectively)
under the relation $\y\x - \x\y = \x$.
Our construction has some
similarity with \cite{C}, where deformation is used to
study the space of harmonic polynomials in two group of
variables. But in \cite{C} after deformation the space of harmonic polynomials
becomes a representation of the double affine Hecke algebra.

The generality we consider in this paper is local Weyl modules over $\gl_r \otimes \CC[x,y]$.
They are labeled by vectors $(\xi_1, \dots, \xi_r)$ such that $\xi_i - \xi_{i+1}$ are
non-negative integers. So up to the action of scalars we can say 
that they are indexed by a partition.

In first section we study finite-dimensional modules. For a partition $\xi$ and an integer 
$N$ we construct a cyclic module $\V(\xi,N)$ over 
$\gl_r \otimes \CC\left<\x,\y\right>$ explicitly.
Then we study its degeneration $V(\xi)$ which is a cyclic module over 
$\gl_r \otimes \CC[x,y]$ not depending on $N$. We show that $V(\xi)$ is a quotient of
the corresponding Weyl module. Our conjecture is that it coincides with the Weyl module.
Here we proof this conjecture for $\xi = (n,s,\dots,s)$, in particular for any weight of
$sl_2$.

Then we find the dimension and describe the action of the constants $\gl_r\otimes 1$ 
on $V(\xi)$ in 
terms of $\rho$-parking functions and the Frobenius characteristic map. 
Also we discuss some character formulas related to 
\cite{HHLRU}.

At last we derive that dimension of the space of diagonal coinvariants is not less than 
the number of parking functions and that the representation of the symmetric group on 
diagonal coinvariants contains the representation on parking functions multiplied by
the sign representation.

In second section we proceed to the limit. First we construct the limits 
of deformed modules using the space of semi-infinite forms. This space is well known 
as a construction for the fundamental representations of the algebra
$\widehat{\gl_r}$. Then we degenerate these inductive systems
and obtain the limits of Weyl modules as representations of a certain central extension of
${\mathfrak s}l_r\otimes \CC[x,x^{-1},y]$. 
It provide us an action of the universal central extension. As a corollary we propose an
explicit formula for the limits of polynomials reciprocal to characters of Weyl modules.

{\bf Acknowledgments}. We are grateful to V.V.Dotsenko and A.N.Kirillov for useful 
and stimulating discussions.
Also we thanks E.B.Feigin for a technical help that made possible to complete this paper being 
at different parts of the world.

BF is partially supported by the grants
RFBR-02-01-01015, RFBR-01-01-00906 and INTAS-00-00055.

SL is partially supported by the grants RFBR-02-01-01015,
RFBR-01-01-00546 and RM1-2545-MO-03.

\section{Deformation of Weyl modules}

\def \mat {{\rm Mat}}

\subsection{Notation}
Let $\mat_r$ be the $r\times r$ matrix algebra. By $E_{ij} \in \mat_r$ denote the matrix units. 
Let $V_r$ be the vector representation of $\mat_r$. 
By $u_i\in V$ denote the basis vectors, so $E_{ij}u_k = u_i$ for $j=k$ and zero otherwise.
Also fix the usual notation ${\rm Id} = E_{11} + \dots + E_{rr}$.

\def \h {{\mathfrak h}}
\def \gs {{\mathfrak sl}}
\def \ep {{\epsilon}}

We will write $\gl_r$ for $\mat_r$ with the usual commutator considered as a Lie algebra.
Let $\h\subset \gl_r$ be the subalgebra spanned by $E_{ii}$, $i=1\dots r$.
Let $\ep_i$, $i=1\dots r$ be the basis in $\h^*$ dual to $E_{ii}$.
Fix the notation $\tau = \ep_1 + \dots + \ep_r$ for the element of $\h^*$ evaluating the trace.

Let $\gs_r \subset \gl_r$ be the subalgebra of traceless matrices.
Choose a basis $h_i = E_{ii} - E_{i+1,i+1}$, $i=1\dots r-1$ 
in the traceless part of $\h$.
By $\omega_i$ denote the dual basis elements in the dual space.

Also note that for an arbitrary associative algebra $A$ we have the structure of 
Lie algebra on the space
$\gl_r \otimes A$ of $A$-valued matrices.
If $A$ is commutative then the notation $\gs_r \otimes A$ makes sense, otherwise
the subspace of traceless matrices is not closed under the commutator.

Concerning the combinatorial notation, for a vector $\eta = (\eta_1, \dots, \eta_m)$ let
$|\eta| = \eta_1+ \dots + \eta_m$. For a partition $\xi = (\xi_1\ge \dots\ge \xi_m \ge 0)$ let
$\xi^t = (\xi^t_1\ge \xi^t_2 \ge \dots)$ be the transposed partition, where
$\xi^t_i$ is the number of $j$ such that $\xi_j\ge i$. At last for a vector
$\eta = (\eta_1, \dots, \eta_m)$ introduce the partition $(1^{\eta_1} 2^{\eta_2} \dots)$, 
where each  $j$ appears $\eta_j$ times.
 
\subsection{Preliminary on deformation and degeneration}

\def \ve {\varepsilon}
\def \gr {{\rm gr}\,}

Let $A$ be an associative (or Lie, etc.) algebra. We say that an algebra ${\mathcal A}$
is a (one parametric) {\em deformation} of $A$ if there exists
a family $A^\ve$ of multiplications (brackets, etc.) on the space 
$A$ depending on $\ve \in \CC$ (that is a family of maps $A\otimes A \to A$) such that

1) $A^0 \cong A$,

2) $A^{\ve} \cong {\mathcal A}$ for $\ve \ne 0$.

\noindent This notion is very useful for various {\em semi-continuous} arguments. 

\def \cb {{\mathcal B}}

Let $\cb$ be an associative (or Lie, etc.) algebra with an increasing filtration 
$\cb_0 \subset \cb_1 \subset \dots$ such that
$\cb_i \cdot \cb_j \subset \cb_{i+j}$ (resp. 
$[\cb_i, \cb_j] \subset \cb_{i+j}$, etc.). 
Then the adjoint graded space $B = \gr \cb$ inherits the structure
of graded algebra (associative, Lie, etc.). This algebra is called a {\em degeneration} of $\cb$.

\begin{prp}
The algebra $\cb$ is a deformation of its degeneration $B$.
\end{prp} 

\begin{proof}
Fix an isomorphism of vector spaces $\phi: B \stackrel{\sim}{\to} \cb$ preserving the filtration.
Let $\deg^\ve : B \to B$ be the map multiplying each graded component $B_i$ by $\ve^i$. Then 
we take $B^0 = B$ and
pullback the operation from the space $\cb$ to
$B^\ve = B$ using the isomorphism $\deg^\ve \circ \phi$ for $\ve \ne 0$. 
\end{proof}

\def \cm {{\mathcal M}}

Let $\cm$ be a module over $\cb$. Suppose $\cm$ is generated by a vector $v$, 
then such a vector is called {\em cyclic}. Introduce the filtration on $\cm$
by setting $\cm_i = \cb_i \cdot v$ (using the universal enveloping algebra for the Lie case). 
The corresponding adjoint graded space $\gr_{v} \cm$
inherits the action of the algebra $B$. We also call this space a {\em degeneration} of $\cm$.

\subsection{Weyl modules}
Let $\g_r = \gl_r\otimes \CC[x,y]$ be the Lie algebra of $r\times r$ matrices over
the polynomial ring in two variables. Clearly 
it splits into the direct sum of the ideal $\gs_r \otimes \CC[x,y]$ of 
traceless matrices and the center $\rm{Id} \otimes \CC[x,y]$.
The elements $E_{ij}^{kl} = E_{ij}\otimes x^k y^l$ for
$1\le i,j \le r$ and  $k,l \ge 0$ form a basis of $\g_r$.

For a vector $\xi = (\xi_1, \dots, \xi_r)$  
define the Weyl module $W(\xi)$ as the maximal finite-dimensional module generated 
by a vector $v_\xi$ such that for any $P \in \CC[x,y]$ we have
\begin{equation}\label{condh}
(E_{ii}\otimes P) \cdot v_\xi = \xi_i P(0,0) \cdot v_\xi, \qquad 
(E_{ij}\otimes P) \cdot v_\xi = 0 \quad \mbox{for} \ 1 \le i < j \le r.
\end{equation}
We can treat $\xi$ as the element $\sum \xi_i \ep_i \in \h^*$, so
$h\otimes P$ acts on $v_\xi$ by the scalar $\xi(h)P(0,0)$ for any $P\in \CC[x,y]$, $h \in \h$. 
It is not immediately clear that such a module exists, so let us deduce it from the
results of \cite{FL2}.

\begin{prp}\label{idp}
For a given $\xi$ consider its {\em signature} 
$\lambda = (\xi_1- \xi_2, \dots, \xi_{r-1}-\xi_r)$.
\begin{enumerate}
\item The module $W(\xi)$ exists.
\item We have $W(\xi) \ne 0$ only if $\lambda$ is dominant (that is all $\lambda_i$
are non-negative integers).
\item For a dominant $\lambda$ the module $W(\xi)$ is isomorphic to 
the module $W_{\CC^2}(\{0\}_\lambda)$ defined 
in \cite{FL2} (while restricted to  $\gs_r \otimes \CC[x,y]$).
\end{enumerate}
\end{prp}

\begin{proof}
The elements ${\rm Id}\otimes P$ belong to the center and therefore 
act on any module generated by $v_\xi$ by the scalars $P(0,0)\cdot (\xi_1 + \dots + \xi_r)$.
So the module $W(\xi)$ can be completely determined by the action of traceless matrices. 

Now take the module $W_{\CC^2}(\{0\}_\lambda)$ with the action of scalars as above. 
In \cite{FL2} it is proved that this module is finite-dimensional,
and it follows from the definition that it is maximal among modules 
generated by $v_\xi$. So we have (iii). 

The statement (ii) follows because otherwise $v_\xi$ can not generate a non-trivial 
finite-dimensional $\gs_r$-module. The statement (i) follows from
the existence of Weyl modules for $\gs_r\otimes \CC[x,y]$.
\end{proof}

For any module $M$ over $\g_r$ and a vector $\eta = (\eta_1, \dots, \eta_r)$ introduce 
the weight subspace $M^\eta$ where $E_{ii}\otimes 1$ acts by $\eta_i$.

\begin{prp} We have
$$W(\xi) = \bigoplus_{{\eta - \xi \in \ZZ^r}\atop{|\eta| = |\xi|}} W(\xi)^\eta,$$
where $|\ |$ returns the sum of coordinates.
\end{prp}

\begin{proof}
Such a statement is always true for a 
finite-dimensional representation of $\gl_r$, where the identity
matrix acts by a scalar. Here we have the action of $\gl_r \otimes 1 \subset \g_r$.
\end{proof}

\def \dh {{\rm D \bf H}}

Let us recall some results from \cite{FL2} about the case 
$\xi = n\ep_1 + s\tau = (n+s,s, \dots,s)$, that is $\lambda = n \omega_1 = (n,0,\dots,0)$.
Following \cite{H} by $\dh_n$ denote the space of diagonal coinvariants
$$\dh_n = \CC[x_1, \dots, x_n; y_1, \dots, y_n]/
\CC[x_1, \dots, x_n; y_1, \dots, y_n]^{\Sigma_n},$$
where the symmetric group $\Sigma_n$ permutes 
the variables $x_i$ and $y_i$ at the same time. The space $\dh_n$ is a representation of 
$\Sigma_n$ as well.

\begin{thm}\label{wsym} \cite{FL2} For an integer $n$, any $s$ and $\eta \in \ZZ^r$
we have an isomorphism of vector spaces
$$W(n\ep_1+s\tau)^{\eta+s\tau} \cong 
\left\{
\begin{array}{ll}
\dh_n^{\Sigma_{\eta_1}\times \dots \times \Sigma_{\eta_r}} & |\eta|=n,\ \eta_i \ge 0;\\
0  & \mbox{otherwise}.
\end{array}\right.$$
\end{thm}

Using the main result of \cite{H} one can calculate  dimensions of these spaces.

\begin{thm}\label{dimc} \cite{FL2}
Dimension of the space $W(n\ep_1+s\tau)^{\eta+s\tau}$ is equal to 
the number of subsets $H \subset \bnr{n}{r}$ such that
$$|H \cap \{1, \dots, n\}\times i| = \eta_i, \quad |H \cap \bnr{k}{r}| \ge k\ \ \mbox{for} \ 
k=1\dots n.$$
\end{thm}

\begin{crl}
We have
$$\dim W(n\ep_1+s\tau) = \frac{(r(n+1))!}{(n+1)!((r-1)(n+1)+1)!},$$
that is the higher Catalan number, and that is
the number of subsets $H \subset \bnr{n}{r}$ such that
$$|H| = n, \quad |H \cap \bnr{k}{r}| \ge k\ \ \mbox{for} \ k=1\dots n.$$
\end{crl}

\subsection{Deformed Weyl modules}

\def \de {\partial}
\def \dif {{\rm Diff}}
\def \cxy {{\CC\left< \x, \y\right>}}

Now let $\cxy$ be the associative algebra generated by the elements $\x$ and $\y$ 
under the relation $\y\x - \x\y = \x$.
It can be considered as the subalgebra of the algebra of differential 
operators $\dif(x)$ in the variable $x$ generated by $\x=x$ and 
$\y =  x \de/\de x$. This subalgebra consists of operators preserving the order of zero 
at the origin.

So $\cxy$ can be degenerated into $\CC[x,y]$ in the standard way.
 Namely let $F^n \cxy$ be the subspace spanned by operators of degree not exceeding $n$,
that is by $\x^i \y^j$ with $j\le n$. 
Then $\gr F^\bullet \cxy \cong \CC[x,y]$.

\def \ggr {{\mathcal G}_r}

Let $\ggr = \gl_r\otimes \cxy$ be the Lie algebra of matrices over $\cxy$ with the
usual commutator. 
In a similar way it can be degenerated into $\g_r$.
For a partition $\xi = (\xi_1 \ge \dots \ge \xi_r \ge 0)$ 
let us define some modules explicitly in order to obtain 
Weyl modules after degeneration.

For an integer $N$ consider the vector space
$V_r^{(N)} = V_r \otimes \left(\CC[x,x^{-1}]/x^N \CC[x]\right)$. 
The vectors $u^i_j = u_j \otimes x^i$ 
for $i < N$, $j=1\dots r$ form a basis of  $V_r^{(N)}$.

The action of $\cxy\subset \dif(x)$ on $\CC[x,x^{-1}]$ preserves the subspace $x^N \CC[x]$,
so the algebra $\cxy$ acts on $\CC[x,x^{-1}]/x^N \CC[x]$.
It provides us the action of $\ggr$ on $V_r^{(N)}$. 
Now introduce the vector 
$$v_\xi^{(N)} = \bigwedge_{j=1}^{r} \bigwedge_{i=N-\xi_j}^{N-1} u_j^i \in 
\bigwedge\nolimits^{|\xi|} V_r^{(N)}.$$
Then define $\V(\xi,N)$ as the $\ggr$-submodule of 
$\bigwedge^{|\xi|} V_r^{(N)}$ generated by the vector $v_\xi^{(N)}$. Note that if $N$ is 
sufficiently large (namely $N\ge \xi_1$) then we 
can avoid negative degrees of $x$.

More generally for a subset $H\subset \NN \times \{1, \dots, r\}$ 
introduce the vector
$$u_H^{(N)} = \bigwedge_{i\times j \in H} u_j^{N-i} \in \bigwedge\nolimits^{|H|} V_r^{(N)}.$$

\begin{prp}\label{condp}
Vectors $u_H^{(N)}$ for 
\begin{equation}\label{condxi}
|H|=|\xi|, \qquad |H \cap \bnr{k}{r}| \ge \xi^t_1 + \dots + \xi^t_k \ \ 
\mbox{for}\  k=1,2,\dots
\end{equation}
form a basis in $\V(\xi,N)$. 
Here $\xi^t = (\xi^t_1 \ge \xi^t_2 \ge \dots)$ denotes the transposed partition.
\end{prp}

\begin{proof}
For $M<N$ introduce the $\ggr$-submodule
$$V_r^{(M\dots N)} = V_r \otimes \left(x^M \CC[x] / x^N \CC[x]\right) \subset V_r^{(N)}.$$ 
We have $\V(\xi,N) \subset \bigwedge^{|\xi|} V_r^{(M\dots N)}$ if $M\le N-\xi_1$.

The action of $\ggr$ on $V_r^{(M\dots N)}$ defines the map  of  algebras $U(\ggr) \to 
\mat_r \otimes \mat_{N-M}$, where $\mat_r$ acts on the factor $V_r$ and $\mat_{N-M}$ acts
on the factor $\left(x^M \CC[x] / x^N \CC[x]\right)$.
 For a pair of integers $M\le i,j < N$ by $F_{ij} \in \mat_{N-M}$
 denote the matrix unit sending $x^i$ to $x^j$.

Let us show that the image of this
map is $\mat_r \otimes B_{N-M}$, where $B_{N-M}$ is the subalgebra of upper-triangular matrices
(spanned by $F_{ij}$ for $i\le j$). 

An element $E_{ij}\otimes \x^k\y^l \in \ggr$ is mapped to
\begin{equation}\label{vnd}
E_{ij} \otimes \left(\sum\limits_{s=M}^{N-k-1}  s^l F_{s,s+k}\right) 
\in \mat_r \otimes B_{N-M}. 
\end{equation}
Therefore any element of $U(\ggr)$ is mapped inside $\mat_r \otimes B_{N-M}$.

Taking different
$l$-s and inverting the Vandermonde matrix
we obtain each $E_{ij} \otimes F_{s,s+k}$ separately from \eqref{vnd} as a linear combination. 
And these elements span
the subalgebra  $\mat_r \otimes B_{N-M}$.

The conditions~\eqref{condxi} mean that $u_H^{(N)} \in \bigwedge^{|\xi|} V_r^{(N)}$ 
and for any $m$
the number of factors $u^i_j$ with $i\ge m$ is not less than in $v^{(N)}_\xi$.
So it is clear that the image of $v_\xi^{(N)}$ under the action of $\ggr$ belongs 
to the linear span of $u_H^{(N)}$ with $H$ satisfying~\eqref{condxi}. 
And any such a monomial can be obtained from $v_\xi^{(N)}$
using the action of matrix elements $E_{ij} \otimes F_{s,s+k}$.
\end{proof}

Now let us degenerate $\V(\xi,N)$ into a $\g_r$-module. 
Note that the degree filtration $F^\bullet \ggr$ 
gives us the corresponding filtration $U^{\le \bullet}(\ggr)$ on the universal enveloping algebra.
Let us choose $v_\xi^{(N)} \in \V(\xi,N)$ as a 
cyclic vector. According to the general procedure  
we project the filtration from $U(\ggr)$ to the cyclic module
$\V(\xi,N)$
and then produce the adjoint graded space.
More precisely 
we set $F^i \V(\xi,N) = U^{\le i}(\ggr) \cdot v_\xi^{(N)}$ and take
$\gr_{v_\xi^{(N)}} \V(\xi, N) = \gr F^\bullet \V(\xi,N)$.

\begin{prp}
The $\g_r$-module $\gr_{v_\xi^{(N)}} \V(\xi, N)$ does not depend on $N$.
\end{prp}

\begin{proof}
Let us construct an isomorphism between the $\g_r$-modules $\gr_{v_\xi^{(N)}} \V(\xi, N)$ and
$\gr_{v_\xi^{(M)}} \V(\xi, M)$.

Let $\phi_N^M: \V(\xi,N) \to \V(\xi, M)$ be the isomorphism mapping $u_H^{(N)}$ to
$u_H^{(M)}$. 
For any $g \in \gl_r$, a set $H$ and integers $i$, $j$, $M$ we have 
$$\phi_N^M\left((g\otimes \x^i \y^j) \cdot u_H^{(N)}\right)
= (g\otimes \x^i (\y+N-M)^j) \cdot u_H^{(M)}.$$
As $\phi_N^M$ maps $v^{(N)}_\xi$ to $v^{(M)}_\xi$ and
the elements $\x^i \y^j$ and $\x^i (\y+N-M)^j$ have the same degree, 
for any $m$ we have 
$$\phi_N^M \left(U^{\le m}(\ggr) \cdot v_\xi^{(N)}\right) 
\subset U^{\le m}(\ggr) \cdot v_\xi^{(M)}.$$
So the map $\phi_N^M$ is compatible with the filtrations, that is 
$\phi_N^M\left(F^m \V(\xi,N)\right) \subset F^m \V(\xi, M)$.

Note that the difference between $\x^i \y^j$ and 
$\x^i (\y+N-M)^j$ is zero after degeneration. Therefore 
$\gr \phi_N^M$ is a map between the adjoint graded modules.

At last note that the map $\gr \phi_M^N$ is inverse to $\gr \phi_N^M$, so we have an 
isomorphism.
\end{proof}

Denote the $\g_r$-module $\gr_{v_\xi^{(N)}} \V(\xi, N)$ for an arbitrary $N$ by $V(\xi)$.

\begin{prp}\label{surj}
There is the canonical surjection $W(\xi) \to V(\xi)$ sending the cyclic vector 
to the cyclic vector.
\end{prp}

\begin{proof}
By definition of Weyl module it is enough 
to show that the cyclic vector of $V(\xi)$ satisfies the conditions \eqref{condh}.

We have $(E_{ij} \otimes \x^m \y^n) \cdot v^{(N)}_\xi =0$ for $i<j$, so the action of 
$E_{ij} \otimes P$ on the cyclic vector of $V(\xi)$ is zero if $i<j$.

The vector $(E_{ii} \otimes \x^m \y^n) \cdot v^{(N)}_\xi$ is
proportional to $v^{(N)}_\xi$ and zero if $m>0$,  so the action of 
$E_{ii} \otimes P$ on the cyclic vector of $V(\xi)$ is zero if $P(0,0)=0$.

At last we have $(E_{ii} \otimes 1) \cdot v^{(N)}_\xi = \xi_i \cdot v^{(N)}_\xi$, so the action of
$E_{ii} \otimes P$ on the cyclic vector of $V(\xi)$ is exactly as required by \eqref{condh}.
\end{proof}

\begin{cnj}\label{wc}
We have $W(\xi) \cong V(\xi)$.
\end{cnj}

\begin{thm}\label{wct}
Conjecture~\ref{wc} is true for $\xi = n\ep_1+s\tau$. 
\end{thm}

\begin{proof}
From Theorem~\ref{dimc} and Proposition~\ref{condp}
we have that $\dim W(n\ep_1+s\tau) = \dim V(n\ep_1+s\tau)$, so 
the surjection is an isomorphism. 
\end{proof}

\subsection{Relation with $\rho$-parking functions}

\def \pf {{\rm PF}}

Let $\rho = (\rho_1 \ge \rho_2 \ge \dots \ge \rho_k \ge 0)$ be a partition. 
A function $f: \{1, \dots n\} \to \{1,\dots n\}$ is called {\em 
$\rho$-parking function} (see \cite{PP}, \cite{PiSt}, \cite{Yan}) if
$$|f^{-1}(\{1, \dots, \rho_{k-s+1}\})| \ge s \ \  \mbox{for}\ s=1\dots k.$$
In particular for $\rho = (n,n-1, \dots,3,2,1)$ we have the usual parking functions (see \cite{H}).

\def \cpf {\CC\pf}

Let $\pf_n(\rho)$ be the set of $\rho$-parking functions.
The symmetric group $\Sigma_n$ acts on the set $PF_n(\rho)$ by
permutation
on the domain. By $\cpf_n(\rho)$ 
denote the corresponding complex representation of 
$\Sigma_n$. Now let us repeat the arguments from \cite{FL2}.

\def \ca {{\mathcal A}}

To describe $\cpf_n(\rho)$ more explicitly introduce the class of
integer sequences
$$\ca_n(\rho) = \{ a_1, \dots, a_n | \ a_i \ge 0, \ a_1+ \dots + a_n = n,\ a_1 +
\dots + a_{\rho_{k-s+1}}
\ge s\ \ \mbox{for}\ s=1\dots k\}.$$

\def \ind {{\rm Ind}}

\begin{prp}\label{indr} We have an isomorphism of $\Sigma_n$-modules
$$\cpf_n(\rho) \cong \bigoplus_{\ca_n(\rho)} \ind_{\Sigma_{a_1} \times \dots \times
\Sigma_{a_n}}^{\Sigma_n} \CC.$$
\end{prp}

\begin{proof}
For a sequence $(a_1, \dots, a_n) \in \ca_n(\rho)$ we take a subset of parking
functions such that $|f^{-1}(\{i\})| = a_i$ for $i=1\dots n$. This set is stable with
respect to the action of $\Sigma_n$ and it forms the permutation
representation isomorphic to
$$\CC[\Sigma_n / \Sigma_{a_1} \times \dots \times
\Sigma_{a_n}] \cong  \ind_{\Sigma_{a_1} \times \dots \times
\Sigma_{a_n}}^{\Sigma_n} \CC.$$
\end{proof}

\def \sgn {{\rm Sign}}

The action of constants $\gl_r \otimes 1$ defines the structure of $\gl_r$-module on
$V(\xi)$. Let us describe it.

\begin{prp}\label{alt} 
For a partition $\xi$ take  $n = |\xi|$ and $\rho = (1^{\xi_1^t}2^{\xi_2^t} \dots )$, 
that is the partition where each $j$ appears $\xi_j^t$ times, 
and $\xi^t$ is the transposed partition. Then we have an isomorphism 
of $\gl_r$-modules 
$$ V(\xi) \cong
\left(V_r^{\otimes n}\otimes \cpf_{n}(\rho)\otimes \sgn\right)^{\Sigma_{n}},$$
where ${\sgn}$ is the sign representation of $\Sigma_n$. At the right hand side $\gl_r$ acts on
$V_r^{\otimes n}$.
\end{prp}

\begin{proof} 
As the filtration on $\V(\xi,N)$ is $\gl_r$-equivariant,  
we have an isomorphism of $\gl_r$-modules  $ V(\xi) \cong \V(\xi,N)$. 
From Proposition~\ref{condp} it follows that
$$\V(\xi,N) \cong \bigoplus_{\ca_n(\rho)} \bigotimes_{i=1}^n \left(\bigwedge\nolimits^{a_i}
V_r\right) \cong \bigoplus_{\ca_n(\rho)}
\hom_{\Sigma_{a_1} \times \dots \times
\Sigma_{a_n}}\left(\sgn
\otimes (V_r^*)^{\otimes n}, \CC \right)$$
as $\gl_r$-modules.
Applying the Frobenius duality and Proposition~\ref{indr} we have
\begin{eqnarray*}
 \V(\xi,N) \cong 
\bigoplus_{\ca_n(\rho)}
\hom_{\Sigma_n}
\left(\sgn
\otimes (V_r^*)^{\otimes n},\ind_{\Sigma_{a_1} \times \dots \times
\Sigma_{a_n}}^{\Sigma_n} \CC \right) \cong  \\
\cong \hom_{\Sigma_n}
\left(\sgn
\otimes (V_r^*)^{\otimes n}, \cpf_n(\rho)\right) \cong
\left(V_r^{\otimes
n}\otimes \cpf_n(\rho)
\otimes \sgn\right)^{\Sigma_n}.
\end{eqnarray*}

\end{proof}

\begin{crl}\label{wtm}
In the notation of Proposition~\ref{alt} for $\eta \in \ZZ^r$ we have
$$\dim V(\xi)^\eta = 
\left\{
\begin{array}{ll}
\dim \left(\cpf_{n}(\rho)\otimes \sgn\right)^{\Sigma_{\eta_1} \times \dots \times \Sigma_{\eta_r}}
& |\eta|=n,\ \eta_i \ge 0;\\
0  & \mbox{otherwise}.
\end{array}\right.
$$
\end{crl}

\begin{proof}
Let us decompose $V_r$ into one-dimensional graded subspaces:
$$\left(V_r^{\otimes n}\otimes \cpf_n(\rho)\otimes \sgn\right)^{\Sigma_n} = 
\bigoplus_{{|\eta|=n}\atop{\eta_i \ge 0}} \left(u_1^{\eta_1} \otimes 
\dots \otimes u_r^{\eta_r}\right)
\otimes  (\cpf_{n}(\rho)\otimes \sgn)^{\Sigma_{\eta_1} \times 
\dots \times \Sigma_{\eta_r}}.$$
\end{proof}

\def \nbl {{\mathcal F}_n^r}

Also we can formulate it in the following way. For a partition $\xi= (\xi_1\ge \dots\ge \xi_m)$
let $\pi_\xi$ be the corresponding
irreducible representations 
of $\Sigma_{|\xi|}$.
And let $\pi^\xi$ be 
the irreducible representation of $\gl_m$ with highest weight $\xi$.

Define the Frobenius map $\nbl$ from
the Grothendiek ring of $\Sigma_n$ to the Grothendiek ring of $\gl_r$ sending $\pi_\xi$
to $\pi^\xi$ if $\xi_{r+1} =0$ and to zero otherwise.

\begin{crl}\label{nabl}
In the notation of Proposition~\ref{alt} we have
an isomorphism of $\gl_r$-modules.
$$V(\xi) \cong \nbl (\cpf_{n}(\rho)\otimes \sgn).$$
\end{crl}

\begin{proof}
The Schur-Weyl duality provide us the decomposition of $\gl_r \times \Sigma_n$-module
$$V_r^{\otimes n} \cong \bigoplus_{{\xi=(\xi_1\ge \dots\ge \xi_r)}\atop{|\xi|=n}}
 \pi_\xi \otimes \pi^\xi.$$
As the irreducible representations of $\Sigma_n$ are real and therefore self-dual, we have
$$(\pi_\xi \otimes \pi_\eta)^{\Sigma_n} = \left\{
\begin{array}{ll}
\CC & \xi = \eta\\
0 & \xi \ne \eta
\end{array}\right.$$
So the statement can be deduced from Proposition~\ref{alt}.
\end{proof}

\subsection{Characters}

\def \ch {{\rm ch}}

The algebra $\g_r$ has the additional 
$\gl_r$-equivariant gradings by degrees of $x$ and $y$. Let us show that the modules $V(\xi)$
are bi-graded.
 
Note that the algebra $\cxy$ and therefore $\ggr$ is
also graded  by degree of $\x$ (but not $\y$). First let us
show that the modules $\V(\xi,N)$ are graded.

\begin{prp}
The module $\V(\xi,N)$ and the subspaces $F^j \V(\xi,N)$ are graded.
\end{prp}

\begin{proof}
For any $H\subset \NN\times \{1,\dots, r\}$ introduce the number
$d(H) = \sum\limits_i i \xi_i^t - \sum\limits_{ i\times j \in H} i$. Then take the linear span
$$\V^m (\xi,N) = \left< \left. u_H^{(N)} \in \V(\xi,N) \right| d(H) = m\right>.$$
We have $\V(\xi,N) = \bigoplus \V^{m}(\xi,N)$ and $g\otimes \x^k \y^l$ maps from
$\V^{m}(\xi,N)$ to $\V^{m+k}(\xi,N)$ for any $g \in \gl_r$.
As the cyclic vector $v^{(N)}_\xi$ belongs to $\V^0 (\xi,N)$,
this is the grading on our $\ggr$-module.

Concerning the filtration, note that the spaces 
$F^m \ggr$ and therefore $U^{\le m} (\ggr)$ are graded.
Therefore we can just take $ F^j \V^m(\xi,N) =  \V^m(\xi,N)\bigcap F^j \V(\xi,N)$ as the graded
components.
\end{proof}

\begin{prp}
The module $V(\xi)$ is bi-graded.
\end{prp}

\begin{proof}
For an arbitrary $N$ let us take 
$$V^{ij}(\xi) = F^j \V^i(\xi,N)/ F^{j-1} \V^i (\xi,N).$$
Then we have $V(\xi) = \bigoplus V^{ij}(\xi)$ 
 and $g\otimes \x^k \y^l$ maps from
$V^{ij}(\xi)$ to $V^{i+k,j+l}(\xi)$ for any $g \in \gl_r$. As the cyclic vector belongs to
$V^{00}(\xi)$, this is the bi-grading on our $\g_r$-module.
\end{proof}

As these gradings are $\gl_r$-equivariant,
we can introduce two additional variables into the character. Let
$$\ch_{x,y} V(\xi)  = \sum_{i,j} x^i y^j \cdot \ch\, V^{ij}(\xi).$$

Now let us calculate the specialization $\ch_x$ of this character when $y =1$.
For any $\rho$-parking function $f$ introduce the integer 
\begin{equation}\label{gre}
|f| =  |\rho| - \sum_{i=1}^n f(i).
\end{equation}
It defines a grading on the space $\cpf_n(\rho)$. So 
$\cpf_n(\rho) = \bigoplus_m \cpf_n^m(\rho)$.

\begin{thm} \label{chmain}
For a partition $\xi$ take  $n = |\xi|$ and $\rho = (1^{\xi_1^t}2^{\xi_2^t} \dots )$, 
that is the partition where each $j$ appears $\xi_j^t$ times, 
and $\xi^t$ is the transposed partition. Then we have 
$$\ch_x V(\xi) = \sum_m x^m \cdot \ch 
\left(V_r^{\otimes n}\otimes \cpf_{n}^m(\rho)\otimes \sgn\right)^{\Sigma_{n}},$$
or equivalently
$$\ch_x V(\xi) = \sum_m x^m \cdot \ch \, \nbl(\cpf_{n}^m(\rho)\otimes \sgn).$$
\end{thm} 

\begin{proof}
Note that for any $N$ we have
$$\ch_x V(\xi) = \sum x^m \cdot \ch\, \V^m(\xi,N).$$

Let us split the set $\ca_n(\rho)$ into subsets $\ca_n^m(\rho)$, containing vectors
$(a_1, \dots, a_n)$ with $|\rho| - \sum_i ia_i = m$. Taking 
into account that $|\rho| = \sum_i i\xi_i^t$, 
similar to the proof of Proposition~\ref{alt} we have
$$\V^m(\xi,N) = \bigoplus_{\ca_n^m(\rho)} 
\bigotimes_{i=1}^n \left(\bigwedge\nolimits^{a_i} V_r\right) \cong 
\nbl(\cpf_{n}^m(\rho)\otimes \sgn).$$
So the grading on $\V(\xi,N)$ coincides
with the grading arising from parking functions. The character formula for $V(\xi)$ follows.
\end{proof}

Note that characters of actual Weyl modules have an additional symmetry.

\begin{prp}
The polynomial
$\ch_{x,y} W(\xi)$ is symmetric with respect to permutation of $x$ and $y$
\end{prp}

\begin{proof}
The Lie group $GL_2$ acts on $\CC[x,y]$ by linear transformations of variables, therefore
it acts on $\g_r$ by the corresponding automorphisms. Note that the
conditions~\eqref{condh} are stable with respect to this action. Therefore $GL_2$ acts
on each Weyl module  $W(\xi)$ preserving $v_\xi$. 

Then note that the gradings on $\g_r$ and therefore on $W(\xi)$ are just 
the action of the diagonal matrix units of $\gl_2 = {\rm Lie}(GL_2)$. So 
the polynomial $\ch_{x,y} W(\xi)$ is character of a $\gl_2$-module, that is a
symmetric polynomial.
\end{proof}

Conjecture~\ref{wc} implies that $\ch_{x,y} V(\xi) = \ch_{x,y} W(\xi)$, so these polynomials
should be symmetric. In particular we know it for $\xi = n \ep_1 + s \tau$.

Then note that the grading~\eqref{gre}
on the space of parking functions coincides with the statistics 
introduced in  \cite{HHLRU} associated with the variable $t$.
Concerning $\ch_{x,y}$ we also believe that the pair of statistics proposed in \cite{HHLRU} 
gives the
correct answer for the usual parking functions (that is $\xi = n\ep_1 + s\tau$).
The second statistic (associated with the variable $q$) makes us thinking
that there exists a generalization of the crystal approach for the two-dimensional case.
But for a general $\xi$ a certain modification is still necessary.

\subsection{On M.~Haiman's theorem}

Note that along the way we proved once more that dimension of 
$\dh_n$ is greater or equal to $(n+1)^{n-1}$ (the number of parking functions) and that 
the representation of $\Sigma_n$ on $\dh_n$ contains the representation on parking functions
multiplied by the sign representation.

To obtain the space of diagonal coinvariants it is enough to consider the zero weight
subspace of the module $W(r\ep_1)$.
Namely Theorem~\ref{wsym} implies that 
$$W(r\ep_1)^{\tau} \cong  \dh_r ^{\Sigma_1 \times \dots \times \Sigma_1} \cong \dh_r.$$
Similarly by Corollary~\ref{wtm} we have 
$$\dim V(r\ep_1)^\tau = |\pf_r|.$$
Actually sets $H$ enumerating the monomials in $\V(r\ep_1,N)^\tau$ 
can be considered as graphics
of parking functions. So the surjection provided by Proposition~\ref{surj}
 gives us the inequality of dimensions.
Note that the proofs of these statements (including Theorem~\ref{wsym}) 
are not based on the results of \cite{H}.

\medskip

Now let us proceed to the action of $\Sigma_r$. Note that we have a statement similar to 
Corollary~\ref{nabl} for the Weyl modules.

\begin{prp}\label{frp}
We have an isomorphism of $\gl_r$-modules
$$W(n\ep_1) \cong \nbl (\dh_n).$$
\end{prp}

\begin{proof}
Theorem~\ref{wsym} implies that character of the $\gl_r$-module
$W(n\ep_1)$ is 
equal to the character
$\left(V_r^{\otimes n}\otimes \dh_n \right)^{\Sigma_{n}}$, so we have
$$W(n\ep_1) \cong \left(V_r^{\otimes n}\otimes \dh_n \right)^{\Sigma_{n}}.$$
Then the statement follows from the Schur-Weyl duality.
\end{proof}

Note that for $r\ge n$ the operator $\nbl$ is 
an inclusion. So let us take $r \ge n$. 
Then as we know that $V(n\ep_1)$ is a quotient of $W(n\ep_1)$, 
combining Corollary~\ref{nabl} and Proposition~\ref{frp} we know that
$\cpf_n\otimes \sgn$ is a subrepresentation of $\dh_n$

\section{Limit of Weyl modules}

\def \sisp {{\mathcal L}}
\def \wsl {{\widehat{\gs_r}}}
\def \gli {{\gl_\infty}}

\subsection{Recall on $\gli$}

Now consider $\infty\times \infty$  matrices, that is matrices $(a_{ij})_{i,j \in \ZZ}$. 
Following \cite{FFu} we say that a matrix $(a_{ij})_{i,j \in \ZZ}$ is a {\em generalized Jacoby} 
matrix if it has only finite number of non-zero diagonals (that is there exists an integer $N$
such that $a_{ij}=0$ when $|i-j|>N$).

Note that we can multiply generalized Jacoby matrices as usual matrices:
$$(a_{ij})(b_{ij}) = \left(\sum_{k\in \ZZ} a_{ik} b_{kj}\right)_{ij}$$ 
and all the sums here are finite. So
introduce the associative algebra $\mat_\infty$ and the Lie algebra $\gli$ of generalized Jacoby
matrices.

\def \sinf {{\infty/2}}
\def \vsinf {\bigwedge\nolimits^{\sinf} V_\infty}

Similar to usual matrices this algebra have the tautological vector representation 
$V_\infty = \left< u_i \right>_{i \in \ZZ}$. 
Consider the space $\vsinf$ spanned by the vectors
$\bigwedge\limits_{i \in H} u_i$ for sets $H \subset \ZZ$ such that the difference of $H$ and
the subset of positive integers $\NN\subset \ZZ$ is finite.

It is clear how to define the action of an element $(a_{ij}) \in \gli$ on $\vsinf$ if $a_{ii}=0$
for all $i$. But for diagonal elements the straightforward approach leads to infinite sums.
And this is the usual situation when a central extension appears.

Consider the subspaces $V_+ = \left< u_i\right>_{i> 0}$ and 
$V_- = \left< u_i\right>_{i\le 0}$ of $V_\infty$. 
Let $\pi: V_\infty \to V_-$ be the projection along $V_+$. Define the {\em Sato-Tate cocycle} by
$$\omega(g_1, g_2) = {\rm Tr} \left(\rule{0pt}{10pt} 
[\pi(g_1), \pi(g_2)] - \pi([g_1,g_2])\right).$$
This notion is well-defined because trace is evaluated on a finite rank matrix.

\begin{thm}{\em (see \cite{FFu})} \label{st}
The central extension $(\gli)^{\widehat{}} = \gli \oplus \CC K$ 
by the Sato-Tate cocycle acts on $\vsinf$ such that
$K$ acts by identity.
\end{thm}

Note that the isomorphism between $(\gli)^{\widehat{}}$ and $\gli \oplus \CC K$ is not canonical. 
But there is the unique way to choose a section $\gli \hookrightarrow (\gli)^{\widehat{}}$
such that the image of diagonal elements acts on the vector 
$$\bigwedge\nolimits^\infty V_+ = \bigwedge\limits_{i>0} u_i \in \vsinf$$ 
by zero and the image of elements with zeroes on the main diagonal acts in the usual way.

\medskip

\def \wg {\widehat{\gl_r}}
\def \wgd {\widehat{\gl_r[D]}}

For our purposes let us fix the isomorphism  $V_r \otimes \CC[x,x^{-1}] \cong V_\infty$
by identifying $u^i_j = u_j \otimes x^i\in V_r \otimes \CC[x,x^{-1}]$ with $u_{ri+j}\in V_\infty$.
Let $\sisp = \bigwedge\nolimits^{\sinf} (V_r\otimes \CC[x,x^{-1}])$.

Similar to the finite-dimensional situation fix the following notation. For a vector
$\xi = (\xi_1, \dots, \xi_r)$ introduce the semi-infinite monomial
$$v^\sinf_\xi =  \bigwedge_{j=1}^r \bigwedge_{i\ge -\xi_j} u^i_j \in \sisp.$$
More generally for $H \subset \ZZ \times \{1, \dots r\}$ such that 
the difference of $H$ and $\NN \times \{1,\dots,r\}$ is finite
introduce the semi-infinite monomial
$$u^{\sinf}_H = \bigwedge_{i\times j \in H} u^{i}_j \in \sisp.$$ 
By definition $\sisp$ is the space spanned by the monomials $u^{\sinf}_H$.

\medskip

Now let us equip $\sisp$ by actions of Lie algebras.
First note that the algebra $\gl_r\otimes\CC[x,x^{-1}]$ acts on $V_r \otimes \CC[x,x^{-1}]$.

\begin{crl}
The central extension $\wg = \gl_r \otimes \CC[x,x^{-1}] \oplus \CC K$ defined by the cocycle
$$\omega(g_1 \otimes x^i, g_2 \otimes x^j) = \delta_{i+j} \cdot {\rm Tr} (g_1 g_2)\cdot i$$
acts on $\sisp$. The element $K$ acts by identity.
\end{crl}

\begin{proof}
Due to the grading any element of $\gl_r\otimes\CC[x,x^{-1}]$ acts by a generalized Jacoby matrix.
So we have a map $\gl_r\otimes\CC[x,x^{-1}] \to \gli$. Then to obtain the action on
 $\sisp$ we just  pullback the Sato-Tate cocycle.
\end{proof}

A bigger but similar example is the algebra of matrix valued 
differential operators $\gl_r\otimes\CC[x,x^{-1},D]$
acting on  $V_r \otimes \CC[x,x^{-1}]$. Here $D = \partial/ \partial x$.

\begin{crl}
The central extension $\wgd = \gl_r\otimes\CC[x,x^{-1},D] \oplus \CC K$ defined by the cocycle
$$\omega(g_1 \otimes x^i D^a, g_2 \otimes x^j D^b) = 
\delta_{i+j-a-b} \cdot {\rm Tr} (g_1 g_2) \cdot \frac{(-1)^a a!b!}{(a+b+1)!} 
\prod_{k=0}^{a+b} (i-k)$$ 
acts on $\sisp$. The element $K$ acts by identity.
\end{crl}

Let $\sisp_n$ be the subspace where the element ${\rm Id}\otimes 1$ acts by the scalar $n$.
Clearly we have the decomposition of the $\wgd$-module (and therefore of the $\wg$-module)
$$ \sisp = \bigoplus_{n \in \ZZ} \sisp_n.$$
Soon we will see that each $\sisp_n$ is irreducible.

\subsection{Recall on representations of $\wg$ and $\wsl$}

\def \he {{\mathcal H}}
\def \th {\tilde{\h}}

Fix the Cartan decomposition $\wg = \n_- \oplus \th \oplus \n_+$. Namely $\n_-$ is spanned by
$E_{ij}\otimes x^s$ for $s<0$ as well as  $s=0$ and $i>j$;
$\n_+$ is spanned by $E_{ij}\otimes x^s$ for $s>0$ as well as  $s=0$ and $i<j$; at last
$\th$ is spanned by $E_{ii}\otimes 1$ and $K$.

Recall that a representation with a highest weight vector is a representation generated by an
eigenvector of $\th$ 
annihilating by $\n_+$. As usual there exists the unique representation of $\wg$ with a given
highest weight, that is
where $K$ acts by a given scalar $k$ and $\h \otimes 1$ acts by a given weight $\xi$ on the 
highest weight vector.
 Let us denote this representation by $L_{k,\xi}$.

Let $\he= \widehat{\gl_1} = {\rm Id}\otimes\CC[x,x^{-1}] \oplus \CC K$ be the 
infinite-dimensional Heisenberg algebra of rank one.
Then we have $\wg = \wsl + \he$ and the intersection of $\wsl$ and $\he$ is the center
$\CC K$.
So any irreducible representation of $\wg$ is the tensor product of an irreducible 
representation of $\wsl$ and an irreducible representation of $\he$.

Any representation of $\he$ with a highest weight vector is isomorphic to the 
Verma module $M_{k,s} = \CC[t_1, t_2, \dots]$, where $K$ acts by the scalar $k$,
the element ${\rm Id}\otimes x^i$ acts by
$t_i$ for $i>0$, by $i \partial/\partial t_i$ for $i<0$ and by the scalar $s$ for $i=0$. 

The algebra $\wsl$ is a {\em Kac-Moody algebra} (see \cite{K}). Note that if $\xi_1-\xi_r \le k$
then the restriction of $L_{k,\xi}$ to $\wsl$ is {\em integrable}. This class of  representations
 generalizes finite dimensional representations of simple algebras. In particular there is
 the action of the affine Weyl group on the weights and 
the Kac-Weyl character formula (see \cite{K})
for irreducible integrable representations. Also integrable representations 
are completely reducible.

\def \xin {\xi(n)}

\begin{prp}
For $n = sr+t$, $0\le t<r$ we have
\begin{equation}\label{xin}
\sisp_n \cong L_{1,\xin} \qquad \mbox{where}\ \xin = s\tau + 
\sum\limits_{i\le t} \ep_i.
\end{equation} 
Here the highest weight vector is identified with $v^\sinf_{\xin} \in \sisp_n$.
\end{prp}

\begin{proof} {\em (sketch)}\ 
First one can check from the definition that $\sisp_n$ is integrable. 
Then the vector $v^\sinf_{\xin}$ generates the submodule isomorphic to $ L_{1,\xin}$.
One can proof that it is the whole module
by comparing the characters.
\end{proof}

The affine Weyl group of $\wsl$ is isomorphic to
$S_r \ltimes Q$, where $Q = \left< \ep_i - \ep_j\right>$ is the root lattice of $\gs_r$.
Similar to the case of $\gl_r$-modules, the Weyl group
acts not only on weights of irreducible representations
but on representations itself.
Let us describe it for $\sisp_n$. 
The action of the subgroup $S_r$ is clear, it just acts on the factor $V_r$.
For $\eta \in \ZZ^r$ define the operators $T_\eta$ sending $u^\sinf_H$ to 
$u^\sinf_{T_\eta H}$, where
$$T_\eta H =   \{ i\times j | (i+\eta_j)\times j \in H\}.$$

One can identify
$Q$ with the set of $\eta \in \ZZ^r$ such that $|\eta|=0$. Then for $\eta \in Q$
the operator $T_\eta$ maps $\sisp_n$ to $\sisp_n$ for any $n$. 
And the action of $T_\eta$ adds $\eta$ to the weight, same as the 
corresponding element of the Weyl group.

\medskip

For our purposes we need the following lemma. Here one can identify
the root $r\omega_1 \in Q$ with $r\ep_1 - \tau \in \ZZ^r$.

\begin{lmm}\label{primc} {\em (see \cite{FS}, \cite{P})}
Let $L$ be an irreducible integrable $\wsl$--module with a highest weight vector $v$. 
By $v(j) \in L$ denote the vector $T_{jr\omega_1} v \in L$.
For $j \ge 0$ take
$F^j L = U(\gs_r \otimes \CC[t])\cdot v(j)$. Then we have
$F^j L \subset F^{j+1} L$ and $L = \bigcup_j F^j  L$.
\end{lmm}

Also 
we will use the 
following construction (see \cite{K}) for modules over $\wg$ at level one.
Note that the subalgebra $\h \otimes \CC[x,x^{-1}]\oplus \CC K$ is isomorphic 
to the infinite-dimensional
Heisenberg algebra $\he_r$ of rank $r$. For $\xi \in \h^*$ by $M_\xi$ denote the Verma
module over $\he_r$, where $K$ acts by identity and $h \otimes 1$ acts by 
the scalar $\xi(h)$ for $h \in \h$. Similar to the rank one case $M_\xi$ is irreducible.

\begin{prp}\label{vop}
We have an isomorphism of $\he_r$-modules
$$\sisp \cong  \bigoplus_{\xi \in \ZZ^r} M_\xi, \qquad 
\sisp_n \cong \bigoplus_{|\xi|=n} M_\xi.$$
Here the highest vector of $M_\xi$ is identified with the vector $v^\sinf_\xi \in \sisp$.
\end{prp}

\begin{proof}{\em (sketch)} \ 
Each vector $v^\sinf_\xi \in \sisp$ generates the $\he_r$-submodule isomorphic to $M_\xi$, so
 the right hand side is included into the left hand side. One can 
show that they are isomorphic by comparing the characters.
\end{proof}

Actually this is a part of the vertex-operator constructions. The action of other elements can
be written in terms of bosonic vertex operators.

\subsection{Inclusions and the limit of deformed modules}
\def \po {\succ}

Let $Q_+ \subset Q$ be the subset of positive roots, that is linear combinations of
$\ep_i - \ep_{i+1}$ with non-negative coefficients.
Introduce the partial order ``$\po$'' on the partitions such that $\xi \po \xi'$
if $\xi-\xi' \in Q_+$. Or equivalently,  $\xi \po \xi'$ 
when $\xi'$ is a weight of the irreducible $\gl_r$-module with highest weight $\xi$.

\begin{prp}\label{pincl1}
Suppose that $\xi \po \xi'$. Then for any $N$ we have
\begin{equation}\label{eincl1}
\V(\xi',N) \subset \V(\xi,N).
\end{equation}
\end{prp}

\begin{proof}
Note that adding $\ep_i - \ep_{i+1}$ to a weight makes
the set of conditions~\eqref{condxi} weaker. 
So the monomials spanning  
$\V(\xi',N)$ are contained in $\V(\xi,N)$.
\end{proof}

\def \det {{\rm Det}}

Let $\det$ be the one-dimensional representation of $\g_r$ evaluating the constant term of 
the trace. Let's deform its powers.
For a pair of integers $M$ and $N$ define the following one-dimensional representation of $\ggr$:
$$\det(M\dots N) =
\left\{
\begin{array}{ll}
\V((N-M)\tau,N),    & N \ge M;\\
\V((M-N)\tau, M)^*, & N \le M.
\end{array}
\right.
$$
Clearly we have $\det(M\dots N)\otimes \det (N\dots K) \cong \det(M\dots K)$ and 
$\gr \det(M\dots N) \cong \det^{N-M}$.
For a shorter notation let's write $\det(N)$ instead of $\det(0\dots N)$.

\begin{prp}\label{detm}
We have 
$$\V(\xi,M) \otimes \det(M\dots N) \cong \V(\xi+(N-M)\tau,N)$$
for $N\ge M$ or when $\xi+(N-M)\tau$ is still a partition.
\end{prp}

\begin{proof}
These spaces can be identified inside $\bigwedge^{|\xi|+(N-M)r} V_r^{(N)}$.
\end{proof}

So we can extend the definition of $\V(\xi,N)$ (and therefore of $V(\xi)$) for an arbitrary
$\xi \in \ZZ^r$. Namely it is zero unless $\xi_1 \ge \xi_2 \ge \dots \ge \xi_r$ and
equal to $\V(\xi-\xi_r \tau,N-\xi_r) \otimes \det(N-\xi_r\dots N)$ otherwise.

\medskip

Note that we have the natural inclusion $\ggr \hookrightarrow \wgd$ 
sending $\x$ to $x$ and $\y$ to $xD$. So we have the action of $\ggr$ on $\sisp$. 

\begin{prp}\label{inclp}
There is the natural inclusion of $\ggr$-modules
\begin{equation}\label{incl}
\V(\xi,N)\otimes \det(-N) \hookrightarrow \sisp_{|\xi|-Nr}
\end{equation}
sending any vector $u\otimes \det(-N)$ from the left hand side to
$u \wedge \bigwedge\limits_{j=1}^r \bigwedge\limits_{i=N}^{\infty} u^i_j$.
\end{prp}

\begin{proof}
We just need to show that this map is a map of representations. The action of $\ggr$ 
on the vector $\bigwedge\limits_{j=1}^r \bigwedge\limits_{i=0}^{\infty} u^i_j$
is zero. So $\bigwedge\limits_{j=1}^r \bigwedge\limits_{i=N}^{\infty} u^i_j \cong \det(-N)$ as
one-dimensional representations of $\ggr$. And $\ggr$ still 
acts on the first factor as on $\V(\xi,N)$.
\end{proof}

Propositions~\ref{pincl1}~and~\ref{detm} provide us the inclusion
$$\V(n\ep_1,N) \otimes \det(N\dots N+1)\hookrightarrow
\V((n+r)\ep_1,N+1)$$
for any $n\ge 0$.
 So we have the inductive system
\begin{equation}\label{is1}
\V(n\ep_1,0) \hookrightarrow \V((n+r)\ep_1,1)\otimes \det(-1)
 \hookrightarrow \V((n+2r)\ep_1,2)\otimes \det(-2) \hookrightarrow \dots
\end{equation}
By Proposition~\ref{inclp} 
each of these spaces is included in $\sisp_n$ as well. Let us consider the inductive limit.

\begin{prp} For the inductive system~\eqref{is1} we have an isomorphism of $\ggr$-modules
$$\lim_\to \left(\rule{0pt}{10pt}\V((n+Nr)\ep_1,N)\otimes \det(-N)\right) \cong \sisp_n$$
\end{prp}

\begin{proof}
We just need to obtain any monomial $u^\sinf_H$ in the right hand side as 
$u^{(N)}_{H'}  \wedge \bigwedge\limits_{j=1}^r \bigwedge\limits_{i=N}^{\infty} u^i_j$ 
for a sufficiently large $N$ and
$u^{(N)}_{H'} \in  \V((n+Nr)\ep_1,N)$. Note that there exists $M$ such that 
$$\{-M,-M+1, \dots \} \times \{1, \dots r\} 
\supset H \supset \{M,M+1, \dots \} \times \{1, \dots r\}.$$
Then for $N\ge M$ we have $u^\sinf_H = u^{(N)}_{H'}  
\wedge \bigwedge\limits_{j=1}^r \bigwedge\limits_{i=N}^{\infty} u^i_j$. And
for $N\ge \frac{r+1}{r-1}M$ the set $H'$ satisfies~\eqref{condxi}, so
we have $u^{(N)}_{H'} \in  \V((n+Nr)\ep_1,N)$.
\end{proof}

\subsection{The Heisenberg filtration}
Now let us construct a filtration on semi-infinite forms compatible with the filtrations on 
deformed Weyl modules.

Define a grading on $\wg$ by setting $\deg (g\otimes x^i) =0$ for a traceless $g$ or $i\ge 0$
and $\deg ({\rm Id} \otimes x^i) =1$ for $i<0$. 
Then introduce the corresponding filtration $F^i \sisp_n = U^{\le i}(\wg) \cdot 
v^\sinf_{\xin}$ on the cyclic module $\sisp_n \cong L_{1,\xin}$. 

According to it we define the following  filtration on $\wgd$. The space $F^i \wgd$ 
contains all differential operators of degree not exceeding $i$ with non-singular symbol of trace.
It means that for $R \in F^i \wgd$ if we write ${\rm Tr}\, R = \sum P_j(x) D^j$ then $P_i(x)$
is non-singular. So we have
\begin{eqnarray*}
E_{ij} \otimes x^sD^t \in F^t \wgd \ \mbox{for}\ i\ne j, &\qquad& h_i \otimes x^sD^t \in F^t\wgd;\\
{\rm Id} \otimes x^sD^t \in F^t \wgd \ \mbox{for}\ s\ge 0, &\qquad& {\rm Id} \otimes x^sD^t \in 
F^{t+1} \wgd;
\end{eqnarray*}
at last the central element $K$ is contained in $F^0 \wgd$.

\begin{prp}
We have $[F^i \wgd, F^j \wgd] \subset F^{i+j} \wgd$.
\end{prp}

\begin{proof}
Degree of the commutator can not exceed the sum of degrees of the entrees.
And symbol of trace of the commutator is just zero.
\end{proof}

Let us show that this filtration is compatible with
the filtration on $\sisp_n$.
For $g \in \gl_r$ and a differential operator $Q$ introduce the generating function (or {\em current})
$$g(Q;z) = \sum_i z^{-i-1} g\otimes (x^iQ).$$
Let us use the normal ordering notation for a product of such currents. 
Namely for $F(z) = \sum z^{-i-1} F_i$ and
$G(z) = \sum z^{-i-1} G_i$ we set
$$:F(z)G(z): = \sum_{i<0,\, j \in \ZZ} z^{-i-j-2} F_i G_j + 
\sum_{i\ge0,\, j \in \ZZ} z^{-i-j-2} G_j F_i.$$
The advantage of it that coefficients of the normal 
ordered product of currents like $g(Q;z)$ are well-defined operators on $\sisp$.
Also by $\partial$ denote the formal derivative of current by $z$.

\begin{lmm}\label{cur}
The coefficients of 
$E_{jj}(D;z)+\frac{1}{2}:E_{jj}(1;z)^2:-\frac{1}{2}\partial E_{jj}(1;z)$  act on $\sisp$
by zero.
\end{lmm}

\begin{proof}
Let us deduce it from Proposition~\ref{vop}. 
First note that this current commutes with 
$h \otimes x^i$ for all $h \in \h$, $i\in \ZZ$. Namely it is clear for $h = E_{ii}$ with $i \ne j$
and we have
$$[E_{jj}(D;z), E_{jj} \otimes x^m] = mz^{m-1} E_{jj}(1;z) -z^{m-2}\frac{m(m-1)}2 K,$$ 
$$[:E_{jj}(1;z)^2:,\, E_{jj}\otimes x^m] = -2mz^{m-1} E_{jj}(1;z),
\qquad [\partial E_{jj}(1;z), E_{jj}\otimes x^m] = -z^{m-2}m(m-1) K.$$
Then by the Schur lemma
the coefficients  act on each subspace $M_\xi$ by a scalar.

It remains to show that these scalars are zeroes. Then we need to check it only for the
action of the grading zero coefficient
(corresponding to $z^{-2}$) on the vectors $v^\sinf_\xi$. It can be easily done using 
$$(E_{jj} \otimes xD) \cdot v^\sinf_\xi = - \frac{\xi_j (\xi_j+1)}2, \qquad
(E_{jj} \otimes 1) \cdot v^\sinf_\xi = \xi_j, \qquad 
(E_{jj} \otimes x^i) \cdot v^\sinf_\xi = 0 \ \ \mbox{for} \ i>0.$$
\end{proof}

\begin{prp}\label{synp}
Elements $h_j\otimes x^i D$ map  $F^l \sisp_n$ into  $F^{l+1} \sisp_n$.
\end{prp}

\def \eqq {\sim}

\begin{proof}
Let us deduce the statement from 
Lemma~\ref{cur}. The elements $h_j\otimes x^i D$ are coefficients of $h_j(D;z)$, therefore they 
act on $\sisp$ in the same way as the coefficients of
$$\frac{1}{2}\left(:E_{j+1,j+1}(1;z)^2: - :E_{jj}(1;z)^2:\right) 
+ \frac{1}{2}\left(\partial E_{jj}(1;z)-\partial E_{j+1,j+1}(1;z)\right).$$
The second summand has values in $U^{\le 1}(\wg)$.
Let us show  that the same is indeed true for the first summand.
As the elements
$E_{jj} \otimes x^t$ commute with the elements $E_{j+1,j+1} \otimes x^s$, it is equal to
$$\frac{1}{2}:\left(E_{j+1,j+1}(1;z)-E_{jj}(1;z)\right)\cdot 
\left(E_{jj}(1;z)+E_{j+1,j+1}(1;z)\right):=
- \frac{1}{2}:h_j(1;z)\cdot \left(E_{jj}(1;z)+E_{j+1,j+1}(1;z)\right):,$$
that has values in $U^{\le 1}(\wg)$ as well.
So the coefficients map $F^l \sisp_n$ into  $F^{l+1} \sisp_n$.
\end{proof}

\begin{prp}\label{sync}
We have $F^k \wgd \cdot F^l \sisp_n \subset F^{k+l} \sisp_n$.
\end{prp}

\begin{proof}
For the elements belonging to $\wg \subset \wgd$ the statement is clear. These are all
the elements of $F^0 \wgd$ and
${\rm Id} \otimes x^i  \in F^1 \wgd$.

For $h_j \otimes x^i D \in F^1 \wgd$  it is shown as Proposition~\ref{synp}.

The elements $E_{st} \otimes x^iD \in F^1 \wgd$ for $s\ne t$ can be obtained by commuting
$h_{\min(s,t)} \otimes x^iD$ and $E_{st}$, therefore the statement for them follows. 
So we completed the proof for $F^1 \wgd$.

At last any element of $F^k \wgd$ can be obtained as a linear combination of 
commutators between elements of 
$F^{k-1} \wgd$ and $F^1 \wgd$. 
Induction on $k$ completes the proof.
\end{proof}

Introduce the notation $V(n+\sinf) = \gr F^\bullet \sisp_n$. 
We have the action of $\gr F^\bullet \wgd$ on $V(n+\sinf)$. Let us describe this algebra in more detail.

\def \sing {{\rm Sing}}

\begin{prp}\label{coc}
The algebra $\gr F^\bullet \wgd$ is an infinite-dimensional central extension of 
$\gs_r \otimes \CC[x,x^{-1},y]$.
\end{prp}  

\begin{proof}
For $g \in \gl_r$, $Q_1,\,Q_2 \in \cxy$ we have
$[{\rm Id} \otimes Q_1, g \otimes Q_2] = g \otimes [Q_1,Q_2]$. As commuting 
decreases the degree of differential operators,
the image of the elements ${\rm Id} \otimes Q$ 
in the adjoint graded space belongs to the center.

Note that in $\wgd$ we have 
$$[g_1\otimes \x^{i_1} \y^{i_1}, g_2\otimes \x^{i_2} \y^{i_2}] = [g_1, g_2] \otimes
\x^{i_1+i_2} \y^{j_1+j_2} \ \ \mbox{modulo} \ \y^{j_1+j_2 -1},$$
so the quotient of $\gr F^\bullet \wgd$ by the image of ${\rm Id} \otimes \cxy$ and $K$ is isomorphic to
 $\gs_r \otimes \CC[x,x^{-1},y]$.
\end{proof}

Identifying the image of ${\rm Id} \otimes \x^i \y^j$ in $\gr F^\bullet \wgd$
with the monomial $x^i y^j$, we can write  
$$\gr F^\bullet \wgd = \gs_r \otimes \CC[x,x^{-1},y] \oplus \CC[x,x^{-1}, y]\oplus \CC K.$$ 
The cocycle defining this extension
 can be written explicitly as
\begin{equation*}
\omega(g_1\otimes P_1, g_2 \otimes P_2) = {\rm Tr}(g_1 g_2)\cdot \sing\{P_1, P_2\} + 
{\rm Tr}(g_1 g_2)\cdot{\rm Res}_x  (dP_1(x,0)\cdot P_2(x,0)) \cdot K
\end{equation*}
where 
$\{,\}$ is the Poisson bracket on $\CC[x,y]$ 
such that $\{y, x\}=1$ and $\sing$ returns the singular part with respect to the variable $x$. 
The first summand appears because ${\rm Id} \otimes \x^i\y^j \in F^{j+1} \wgd$ for $i<0$, 
the second summand remains from the Sato-Tate cocycle.

\subsection{Back to the Weyl modules}
Now let us proceed from $\ggr$-modules $\V(\xi, N)$ and $\sisp_n$ to the adjoint graded
$\g_r$-modules $V(\xi)$ and $V(n+\sinf)$.

\begin{lmm}\label{ord}
Suppose that we have two partitions $\xi \succ \xi'$ (that is  $\xi-\xi' \in Q_+$). 
Then the vector $v^{(N)}_{\xi'} \in \bigwedge^{|\xi|} V^{(N)}_r$ 
can be obtained from $v^{(N)}_{\xi}$ 
by the action of $\gs_r \otimes \CC[x]$.
\end{lmm}

\begin{proof}
First consider
the case $\xi - \xi' = \ep_i - \ep_j$, $i<j$. Then $\xi_i - \xi_j' > 0$ and we have
\begin{equation}\label{tmp}
v^{(N)}_{\xi'} = \left(E_{ji}\otimes x^{\xi_i - \xi_j'}\right) \cdot v_\xi^{(N)}.
\end{equation}

In general we have  $\xi - \xi' = \sum_s \ep_{i_s} - \ep_{i_s+1}$
for some integers $i_s$. So we just apply~\eqref{tmp} simultaneously for each summand.
\end{proof}

\begin{prp}\label{prf}
The inclusion~\eqref{eincl1} is compatible with the filtrations and thereby induces the map
$V(\xi') \to V(\xi)$.
\end{prp}

\begin{proof}
By Lemma~\ref{ord} the cyclic vector $v^{(N)}_{\xi'}$ is mapped into 
$U(\gs_r \otimes \CC[x]) \cdot v^{(N)}_{\xi} \subset  F^0 \V(\xi,N)$. So any vector
from $F^i \V(\xi',N) = U^{\le i}(\ggr) \cdot v^{(N)}_{\xi'}$ is mapped into
$U^{\le i}(\ggr) \cdot U(\gs_r \otimes \CC[x]) \cdot v^{(N)}_{\xi} \subset F^i \V(\xi,N)$.

Taking the adjoint graded spaces we obtain the adjoint graded map $V(\xi') \to V(\xi)$.
\end{proof}

\begin{crl}
We have the natural
map
$$\alpha_{n,s}: W(n\ep_1+ s\tau) \to W((n+r)\ep_1+(s-1)\tau)$$
sending 
$v_{n\ep_1+ s\tau}$ to 
$$\left(E^{n+1,0}_{r1} \dots E^{n+r-2,0}_{31} E^{n+r-1,0}_{21}\right) \cdot 
v_{(n+r)\ep_1+(s-1)\tau}.$$
\end{crl}

\begin{cnj}
The map $\alpha_{n,s}$ is an inclusion. 
\end{cnj}

Anyway the degeneration of the inductive 
system~\eqref{is1} is
\begin{equation}\label{is2}
W(n\ep_1) \stackrel{\alpha_{n,0}}{\longrightarrow}
W((n+r)\ep_1-\tau) 
\stackrel{\alpha_{n+r,-1}}{\longrightarrow} W((n+2r)\ep_1-2\tau) 
\stackrel{\alpha_{n+2r,-2}}{\longrightarrow} \dots.
\end{equation}

\begin{prp}\label{grincl}
The map~\eqref{incl} is compatible with the filtrations and thereby induces the
map $V(\xi) \to  V(|\xi|+\sinf)$.
\end{prp}

\begin{proof}
Concerning the filtration, by Proposition~\ref{sync} 
it is enough to show that $v^{(N)}_\xi\otimes \det(-N)$ is mapped into $F^0 \sisp_{|\xi|-Nr}$.
The image of this vector is exactly $v^{\sinf}_{\xi - N\tau}$, that belongs to 
the Weyl group orbit of the highest weight vector. Therefore it can be obtained from the highest 
weight vector by the action of the zero graded subalgebra $\wsl\subset \wg$.

After degeneration we have the map $V(\xi - N \tau) \to V(|\xi| - Nr + \sinf)$. 
Note that this map will be the same if we add $M\tau$ to $\xi$ and $M$ to $N$. 
So we just take $N=0$.
\end{proof}

\begin{crl}
There is the natural map $W(n\ep_1+s\tau) \to  V(n+rs+\sinf)$.
\end{crl}

\begin{lmm}\label{ord2}
Suppose that we have two vectors $\xi$ and $\xi'$ such that $\xi-\xi' \in Q_+$ 
Then 
the vector $v^\sinf_{\xi'} \in \sisp$ 
can be obtained from 
$v^\sinf_{\xi}$ 
by the action of $\gs_r \otimes \CC[x]$.
\end{lmm}

\begin{proof}
Similar to Lemma~\ref{ord}
\end{proof}

\begin{prp}\label{grlim}
For the inductive system obtained from~\eqref{is1} and any $i$  we have
\begin{equation}\label{flim}
\lim_\to \left(F^i \V((n+Nr)\ep_1,N)\otimes \det(-N)\right) \cong F^i \sisp_n
\end{equation}
\end{prp}

\begin{proof}
From Proposition~\ref{grincl} it follows that the left hand side  is contained in the
right hand side. Let us show the equality.

1) {\em Case $i=0$}. Let $L=F^0 \sisp_n$. 
It is an irreducible $\wsl$--module with the highest weight vector
$v(0) = v^\sinf_{\xin}$. Let us deduce the statement from Lemma~\ref{primc}.

The image of $v^{(N)}_{(n+Nr)\ep_1}$ in $L$
is equal to $v^\sinf_{(n+Nr)\ep_1-N\tau}$. The vector 
$v(j) = T_{jr\omega_1} v^\sinf_{\xin}$ is equal to
$v^\sinf_{\xin +jr\ep_1-j\tau}$.

Let $n = sr+t$, $0\le t<r$. Then we have 
$$t\ep_1+s\tau -\xin = \sum_{k=2}^t (\ep_1 - \ep_k) \in Q_+,$$
so for $N\ge j-s$ the difference between $(n+Nr)\ep_1-N\tau$ and $\xin +jr\ep_1-j\tau$ belongs to
$Q_+$.

As the action of $\gs_r\otimes \CC[x]$ preserves $F^0$, 
Lemma~\ref{ord2} implies that $v(j)$ and therefore $U(\gs_r\otimes \CC[x]) \cdot v(j)$
belongs to the image of the left hand side.
Tending $j \to \infty$ and applying Lemma~\ref{primc} we obtain the statement~\eqref{flim}
for $F^0$.

2) {\em General Case}.
Let $I_k = {\rm Id} \otimes x^{-k}$. As $I_k$ commutes with $\wsl$ we can 
reformulate the statement in the following way. Let $P$ be a polynomial in variables
$I_k$, $k>0$ of degree $i$. Then we need to show that any element of $P \cdot L$ belongs to
the image of $F^{i} \V((n+Nr)\ep_1,N)$ for a certain $N$. Due to Lemma~\ref{primc}
it is enough to obtain the vectors $P\cdot v(j)$.

Now suppose we know this statement for $F^{i-1}$. Let us show it for $F^i$. Decompose
$P = \sum_k I_k P_k$, such that $\deg P_k \le i-1$. Then we need to obtain any element 
$I_k P_k \cdot v(j)$ from the left hand side of~\eqref{flim} .

\def \htt {{h_+}}

Let $\htt = h_1\otimes x^2D \in F^1 \wgd$.
We have $[\htt, I_s] = -s\cdot h_1\otimes x^{1-s}$ and $\htt \cdot v(j) =0$ because 
its weight (in the sense of $\wsl$) is outside the set of weights for
$\sisp_n$. So $\htt P_k \cdot 
v(j)$ belongs to $F^{i-1} \sisp_n$. Also we have
$$[\htt, h_1\otimes x^{-k-1}] = -(k+1)(E_{11}+E_{22})\otimes x^{-k} =
-\frac{2k+2}{r}I_k + \mbox{lower terms},$$
where lower terms belong to $\wsl = F^0 \wgd$. Therefore
$$\htt P_k (h_1\otimes x^{-k-1}) \cdot  v(j) = 
\htt (h_1\otimes x^{-k-1}) P_k \cdot  v(j) =
-\frac{2k+2}{r} I_k P_k \cdot v(j) + \mbox{lower terms},$$
where lower terms belong to $F^{i-1} \sisp_n$. 

By assumption for each $j$ there exists  $N$ such that $P_k (h_1\otimes x^{-k-1})  \cdot v(j)$ 
is the image of a certain element $u \in F^{i-1} \V((n+Nr)\ep_1,N)$. 
Then the image of 
$-\frac{r}{2k+2} \cdot(h_1\otimes \x^2\y) \cdot u \in F^i \V((n+Nr)\ep_1,N)$ is 
$I_k P \cdot v(j)$ modulo
 $F^{i-1} \sisp_n$. Applying the assumption once more for the lower
terms, we obtain the statement for $F^i$.
\end{proof}

\begin{lmm}\label{lalg}
Suppose we have an inductive system $A_1 \to A_2 \to \dots$ and subspaces $B_i \subset A_i$
such that each $B_i$ is mapped inside $B_{i+1}$. Then we have the canonical isomorphism
$$\frac{\lim\limits_\to A_i}{\lim\limits_\to B_i} \cong \lim_\to \left(A_i / B_i\right).$$
\end{lmm}

\begin{thm} For the inductive system~\eqref{is2} we have an isomorphism of $\g_r$-modules
$$\lim_\to W((n+Nr)\ep_1 -N\tau) \cong  V(n+\sinf).$$
\end{thm}

\begin{proof}
Using Proposition~\ref{grlim} then Lemma~\ref{lalg} and then Theorem~\ref{wct}
we have $V(n+\sinf) = \gr F^\bullet \sisp_n=$
\begin{eqnarray*}
=\bigoplus_i \frac{F^i \sisp_n}{F^{i-1} \sisp_n} \cong 
 \bigoplus_i \frac{\lim\limits_\to \left(F^i \V((n+Nr)\ep_1,N)\otimes\det(-N)\right)}
{\lim\limits_\to \left(F^{i-1} \V((n+Nr)\ep_1,N)\otimes\det(-N)\right)} 
\cong \bigoplus_i \lim_\to \frac{F^i \V((n+Nr)\ep_1,N)\otimes\det(-N)}
{F^{i-1} \V((n+Nr)\ep_1,N)\otimes\det(-N)}  =\\
= \lim_\to \left(V((n+Nr)\ep_1) \otimes \det^{-N}\right) 
\cong \lim_\to W((n+Nr)\ep_1-N\tau),
\end{eqnarray*}
that is the limit of the inductive system~\eqref{is2}.
\end{proof}

To write this statement in a more elegant form one can return to the notation of \cite{FL2}.
Recall that by Proposition~\ref{idp} we have $W((n+Nr)\ep_1-N\tau) \cong 
W_{\CC^2}(\{0\}_{(Nr+n)\omega_1})$ as  $\gs_r\otimes \CC[x,y]$-modules.

\begin{crl} 
The action of $\gs_r\otimes \CC[x,y]$ 
on the limit  of the modules $W_{\CC^2}(\{0\}_{(Nr+n)\omega_1})$
can be extended to an action of the
universal central extension of $\gs_r \otimes \CC[x,x^{-1},y]$.
\end{crl}

\begin{proof}
Taking into account Proposition~\ref{coc} we know that the action of  $\gs_r\otimes \CC[x,y]$ 
can be extended to an action of a certain central extension of $\gs_r \otimes \CC[x,x^{-1},y]$.
Therefore it can be extended to an action of the universal one.
\end{proof}

\begin{crl}
The action of $\gs_r\otimes \CC[x] \subset \gs_r\otimes \CC[x,y]$ on the limit  
of the modules $W_{\CC^2}(\{0\}_{(Nr+n)\omega_1})$ can be extended to an action of
$\wg = \wsl+{\mathcal H}$.
\end{crl}

\begin{proof}
As the Heisenberg filtration is $\wsl$-equivariant, the spaces $\sisp_n$ and
$V(n+\sinf)$ are isomorphic as $\wsl$-modules.
\end{proof}

Note that the isomorphism between $\wsl$-modules
$\sisp_n$ and $V(n+\sinf)$ can be established in the canonical way
by identifying the action of ${\rm Id} \otimes x^{-i}$ on $\sisp_n$ with the action of
the central element corresponding to $x^{-i}$ on $V(n+\sinf)$.

\subsection{Limit of characters}
Another consequence of these considerations is an explicit formula for the limits of characters
of Weyl modules.
First let us count their degrees. Let
$$d(\xi) = \deg_x \ch_{x,y} V(\xi).$$

\begin{prp}\label{degp}
We have 
$$d(\xi) = \sum_i i\xi^t_i - \sum_{i=1}^{|\xi|} \left[ \frac{i+r-1}{r} \right],$$
where $[x]$ returns the integer part of $x$.
\end{prp}

\begin{proof}
Note that $\deg_x \ch_{x,y} V(\xi) = \deg_x \ch_{x} V(\xi)$, so we can use
Theorem~\ref{chmain}. A $\rho$-parking function $f$ contributes to the character only if
 $|f^{-1}(\{j\})|\le r$ for all $j$.
So the $\rho$-parking function $f$ such that
$f(i) = \left[ \frac{i+r-1}{r} \right]$ contributes to the maximal degree.
\end{proof}

In particular by Theorem~\ref{wct} we have
$$ \deg_x \ch_{x,y} W(n\ep_1 + s\tau) = d(n\ep_1+s\tau) =
\sum_{i=1}^n \left(i - \left[ \frac{i+r-1}{r} \right]\right).$$

Note that the 
$\ggr$-module $\sisp$ is also graded by degree of $\x$. Namely for any 
$H \subset \ZZ \times \{1, \dots, r\}$, such that the difference of $H$ and
$\NN \times \{1, \dots, r\}$ is finite,
 introduce the non-negative integer
$$d(H) = 
 \sum_{{i\times j \notin H}\atop{i\ge 0}} i -
\sum_{{i\times j \in H}\atop{i<0}} i.$$
Let $\sisp^m$ be the subspace spanned by $u^\sinf_H$ with $d(H)=m$. Then 
$g \otimes \x^i\y^j$ maps from $\sisp^m$ to $\sisp^{m-i}$ for any $g \in \gl_r$. 

The submodules $\sisp_n$ are graded. Let us take $\sisp_n^m = \sisp_n \bigcap \sisp^{d_n + m}$,
where $d_n$ is degree of the highest weight vector $v^\sinf_{\xi(n)}$ of $\sisp_n$.
So the grading on
  $\sisp_n$ starts from $\sisp_n^0$. 

\begin{prp}\label{invp}
The inclusion~\eqref{incl} maps $\V^i(\xi,N) \otimes \det(-N)$ to 
$\sisp_{|\xi|-Nr}^{d(\xi) -i}$.
\end{prp}

\begin{proof}
Note that the image of the cyclic vector $v^{(N)}_\xi \otimes \det(-N)$ is graded,  
let $d$ be its degree.
Then the subspace $\V^i(\xi,N) \otimes \det(-N)$ is mapped to $\sisp_{|\xi|-Nr}^{d -i}$. 
One can show that $d=d(\xi)$ by an explicit calculation using Proposition~\ref{degp}
or just by observation that the highest weight vector of the right hand side 
belongs to the image of the left hand side.
\end{proof}

\def \rec {{\rm Rec}_x}

Now let us present a formula for the limit of polynomials, reciprocal to characters
of $W(n\ep_1 + s\tau)$ with respect to the variable $x$. Introduce the notation
$$\rec (P(x, y)) = x^{\deg_x P} P(x^{-1}, y).$$

\begin{thm} We have
$$\lim_{N\to \infty} \rec  \left(\rule{0pt}{10pt}\ch_{x,y} W((n+Nr)\ep_1-N\tau)\right) =
\frac{\sum\limits_{{\xi \in \ZZ^r}\atop{|\xi|=n}} {e^\xi} x^\frac{(\xi,\xi)- (\xin, \xin)}2
}{\prod\limits_{i>0}
(1-x^iy)(1-x^i)^{r-1}},$$
where $(\xi, \xi) = \sum_i \xi_i^2$ and $\xin$ is the highest weight of $\sisp_n$ 
(see~\eqref{xin} for the explicit formula).
\end{thm} 

\begin{proof}
First let us show that the right hand side is equal to
\begin{equation}\label{rhsm}
\sum_{i,j} x^i y^j\cdot \left(\ch \, F^j \sisp^i_n - 
\ch \, F^{j-1} \sisp_n^i\right).
\end{equation}
To do it let us use the decomposition provided by
Proposition~\ref{vop}. As the vectors $v^\sinf_\xi$ are graded and belongs to
$F^0$, this decomposition is compatible with the grading and the filtration. Then we have
$$\sum_{i,j} x^i y^j\cdot \left(\ch \, F^j M^i_\xi - 
\ch \, F^{j-1} M^i_\xi\right) = 
\frac{e^\xi x^{\deg(v^\sinf_\xi)}}{\prod\limits_{i>0} (1-x^iy)(1-x^i)^{r-1}},$$
namely factors $(1-x^i)$ correspond to the action of $h_j\otimes x^{-i}$ 
and the factor $(1-x^iy)$ corresponds to the action of ${\rm Id}\otimes x^{-i}$.
Summing these characters up we obtain the right hand side.

Due to Theorem~\ref{wct} for any $N$ we have
$$\ch_{x,y} W(m\ep_1) = \sum_{i,j} x^i y^j\cdot \left(\ch \, F^j \V^i(m\ep_1,N) - 
\ch \, F^{j-1} \V^i(m\ep_1,N)\right),$$
and therefore for any integer $s$ the character $\ch_{x,y} W(m\ep_1+ s\tau)$ is equal to
$$\sum_{i,j} x^i y^j\cdot \left(\ch (F^j \V^i(m\ep_1,N)\otimes \det(N\dots N+s)) - 
\ch (F^{j-1} \V^i(m\ep_1,N)\otimes \det(N\dots N+s)) \right).$$

Then using the 
inclusion~\eqref{incl} (here according to Proposition~\ref{invp} 
we invert the variable $x$ in the character)
and taking into account Proposition~\ref{grlim} we proceed to the limit and
obtain character~\eqref{rhsm} of $\sisp_n$. 
\end{proof}

In particular it can be considered as another test for the statistics proposed in \cite{HHLRU}.


\begin{thebibliography}{}

\bibitem[CL]{CL} V.~Chari, T.~Le, {\em Representations of Double Affine  Lie algebras}, math.QA/0205312.

\bibitem[CP]{CP} V. Chari, A. Pressley, {\em Weyl Modules for Classical 
and Quantum Affine Algebras}, 
Represent. Theory 5 (2001), 191--223, math.QA/0004174.

\bibitem[C]{C} I~Cherednik, {\em Diagonal Coinvariants and Double Affine Hecke Algebras}, to appear in IMRN, 
math.QA/0305245.

\bibitem[FF]{FF}
B.~L. Feigin and E.~Feigin,
{\em $q$-characters of the tensor products in
$\gs_2$-case},
math.QA/0201111.

\bibitem[FFu]{FFu}
B.~L. Feigin and D.~B.~Fuchs,
{\em Representations of the Virasoro algebra.  Representation of Lie groups and related topics},
465--554, Adv. Stud. Contemp. Math., {\bf 7}, Gordon and Breach, New York, 1990.

\bibitem[FL1]{FL1} B. Feigin, S. Loktev, {\em On Generalized Kostka
Polynomials and the Quantum Verlinde Rule}, Differential
topology, infinite--di\-men\-si\-onal Lie algebras, and applications,
Amer. Math. Soc. Transl. Ser. 2, Vol. 194 (1999), p. 61--79,
math.QA/9812093.

\bibitem[FL2]{FL2} B. Feigin, S. Loktev, {\em Multi-dimensional Weyl modules and Symmetric Functions},
to appear in Comm. Math. Phys., math.QA/0212001.


\bibitem[FS]{FS} B. Feigin and A. Stoyanovsky, 
{\it Quasi-particles
models for the representations of Lie algebras and geometry of flag
manifold}, hep-th/9308079, RIMS 942; {\it Functional models for the
representations of current algebras and the semi-infinite Schubert
cells}, Funct. Anal. Appl. {\bf 28} (1994), 55--72.


\bibitem[HHLRU]{HHLRU} J. Haglund, M. Haiman, N. Loehr, J. B. Remmel, A. Ulyanov,
{\em A Combinatorial Formula for the Character of the Diagonal Coinvariants}, math.CO/0310424.

\bibitem[H]{H} M. Haiman, {\em Vanishing theorems and character formulas
for the Hilbert scheme of points in the plane}, Invent. Math. 149, no. 2
(2002) 371-407, math.AG/0201148.

\bibitem[K]{K} V. Kac, {\em Infinite dimensional Lie algebras}, Cambridge
Univ. Press, 1985.


\bibitem[P]{P} M. Primc,
{\it Vertex operator construction of standard modules for 
$A^{(1)}_n$},
Pacific J. Math. {\bf 162} (1994) 143--187;
{\it Standard representations of $A^{(1)}_n$}, 
in Proceedings of Marseilles Conference  
`Infinite Dimensional Lie Algebras and Groups',
Advanced Studies in Mathematical Physics,  
World Scientific, Singapore, 
{\bf 7} (1988), 273--282.



\bibitem[PP]{PP} I.~M.~Pak, A.~E.~Postnikov,
{\em Resolutions for $S_n$-modules corresponding to skew hooks, 
and combinatorial applications},  
Functional Analysis and its Applications {\bf 28} (1994), no.~2, 132--134.

\bibitem[PS]{PiSt} J.~Pitman, R.~Stanley, {\em  A polytope related to 
empirical distributions, plane trees, parking functions, and the associahedron},
Discrete and Computational Geometry {\bf 27} (2002), 603--634. 

\bibitem[S]{S} R. Stanley, {\em Enumerative Combinatorics}, Cambridge
Studies in Advanced Mathematics, Vol. {\bf 62}.


\bibitem[Y]{Yan} C.~H.~Yan,
{\em On the enumeration of generalized parking
functions}, Proceedings of the 31-st Southeastern International
Conference on Combinatorics, Graph Theory and Computing 
(Boca Raton, FL, 2000),
Congressus Numerantium {\bf 147} (2000), 201--209.



\end{thebibliography}
\end{document}